\documentclass[11pt,a4paper] {amsart}
\usepackage{amssymb}
\usepackage{amsthm}
\usepackage{latexsym,mathrsfs}
\usepackage{pifont}

\newtheorem{theorem}{T{\hskip 0pt\footnotesize\bf HEOREM}}[section]
\newtheorem{lemma}[theorem]{L{\hskip 0pt\footnotesize\bf EMMA}}
\newtheorem{proposition}[theorem]{P{\hskip 0pt\footnotesize\bf ROPOSITION}}

\newtheorem{corollary}[theorem]{C{\hskip 0pt\footnotesize\bf OROLLARY}}
\newtheorem{remark}[theorem]{R{\hskip 0pt\footnotesize\bf EMARK}}

\newcommand{\bprop} {\begin{proposition}}
\newcommand{\eprop} {\end{proposition}}
\newcommand{\btheo} {\begin{theorem}}
\newcommand{\etheo} {\end{theorem}}
\newcommand{\blem} {\begin{lemma}}
\newcommand{\elem} {\end{lemma}}
\newcommand{\bcor} {\begin{corollary}}
\newcommand{\ecor} {\end{corollary}}

\newcommand{\Be}{\begin{equation}}
\newcommand{\Ee}{\end{equation}}
\newcommand{\Bea}{\begin{eqnarray}}
\newcommand{\Eea}{\end{eqnarray}}
\newcommand{\Bes}{\begin{equation*}}
\newcommand{\Ees}{\end{equation*}}
\newcommand{\Beas}{\begin{eqnarray*}}
\newcommand{\Eeas}{\end{eqnarray*}}
\newcommand{\Ba}{\begin{array}}
\newcommand{\Ea}{\end{array}}

\def\C{\mathbb{C}}
\def\N{\mathbb{N}}

\def\B{\mathcal B}

\scrollmode

\begin{document}




\title[Logarithmic Bloch spaces in the polydisc]{Logarithmic Bloch spaces in the polydisc, endpoint results for Hankel operators and pointwise multipliers}
 \author{Beno\^it F. Sehba}
 \address{Department of Mathematics, University of Ghana, PO. Box LG 62 Legon, Accra Ghana}

\maketitle


\begin{abstract}
We define two notions of Logarithmic Bloch space in the polydisc for which we provide equivalent definitions in terms of symbols of bounded Hankel operators. We also provide a full characterization of the pointwise multipliers between two different  Bloch spaces of the unit polydisc.
\end{abstract}



\section{Introduction}\label{section1}
For $0<p<\infty$, the Bergman space $A^p(\mathbb{D})$ of the unit disc $\mathbb{D}$ of the complex plane $\mathbb{C}$ consists of all holomorphic functions $f$ on $\mathbb{D}$ such that
$$\|f\|_p:=\left(\int_{\mathbb{D}}|f(z)|^pd\nu(z)\right)^{1/p}<\infty$$
where $\nu$ is the normalized Lebesgue measure on $\mathbb{D}$. The orthogonal projection from $L^2(\mathbb{D})$ onto its closed subspace $A^2(\mathbb{D})$ is called the Bergman projection and denoted $P$.
\medskip

For a bounded holomorphic function on $\mathbb{D}$, the Hankel operator with symbol $b$ is the operator defined for any integrable function $f$ on $\mathbb{D}$ by 
\begin{equation}\label{eq:hakeldef}
    h_b(f):=P(b\overline{f}).
\end{equation}

Recall that a holomorphic function $f$ on $\mathbb{D}$ is said to be a Bloch function if $$\sup_{z\in\mathbb{D}}(1-|z|^2)|f'(z)|<\infty.$$
The set of all Bloch functions is called the Bloch space and denoted $\mathcal{B}(\mathbb{D})$. It is a Banach space when endowed with the following norm
$$\|f\|_{\mathcal{B}}:=|f(0)|+\sup_{z\in\mathbb{D}}(1-|z|^2)|f'(z)|.$$
The Bloch space in one-parameter can be identified as the dual space of the Bergman space $A^1(\mathbb{D})$ (see \cite{Zhu}). Equivalent definitions of this space are also given in terms of symbols of bounded Hankel operators on the Bergman spaces $A^p$, $1<p<\infty$ (see for example \cite{JPR}) and image of $L^\infty(\mathbb{D})$ by the Bergman projection (see \cite{Zhu}). These equivalent characterizations of the Bloch space in one-parameter extend to higher-parameter \cite{Constantin, stevo, Zhu1}. 
\medskip

Also in the unit disc, an analytic function $b$ is a multiplier of the Bloch space if and only if it is bounded and satisfies the following Bloch-type condition 
$$\sup_{z\in\mathbb{D}}(1-|z|^2)\left(\log\left(\frac{2}{1-|z|^2}\right)\right)|f'(z)|<\infty.$$
The space of all holomorphic functions satisfying the last condition is the sometimes called logarithmic Bloch space and denoted $\mathcal{B}_{L}(\mathbb{D})$ (more often $\mathcal{B}_{log}(\mathbb{D})$). We endowed it with the norm
$$\|f\|_{\mathcal{B}_{L}}:=|f(0)|+\sup_{z\in\mathbb{D}}(1-|z|^2)\left(\log\left(\frac{2}{1-|z|^2}\right)\right)|f'(z)|.$$
It is also known that $\mathcal{B}_{L}(\mathbb{D})$ is the exact range of symbols of bounded Hankel operators on $A^1(\mathbb{D})$ (see for example \cite{BL}).
\medskip

In the multi-parameter case, i.e. the polydisc $\mathbb{D}^n$, there is one Bloch space that corresponds to the dual space of the Bergman space $A^1(\mathbb{D}^n)$ defined in the next section, and another Bloch space which is a subspace of the first one and can be obtained using the definition of a Bloch space in several complex variables by R. Timoney \cite{Timoney}. For each of these Bloch spaces, we introduce a corresponding logarithmic Bloch space for which we provide an equivalent definition in terms of set of symbols of bounded Hankel operators.  We also characterize the multipliers algebra of the Bloch space corresponding to the dual of $A^1(\mathbb{D}^n)$ and the pointwise multipliers from the smaller Bloch space to the latter.
\section{Function spaces}\label{section2}
Recall that for $0<p<\infty$, the Bergman space $A(\mathbb D^n)$ consists of analytic functions $f$ in $\mathbb D^n$ such that
\Be
\|f\|_{p}^p:=\int_{\mathbb D^n}|f(z)|^pd\nu(z)<\infty,
\Ee

here $d\nu(z)=d\nu_1(z_1)\cdots d\nu_n(z_n)$ for $z=(z_1,\ldots,z_n)$, where $d\nu_j$ is the normalized Lebesgue measure on the unit disc $\mathbb D$.  That is $A^p(\mathbb D^n)$ is the subspace of the Lebesgue space $L^p=L^p(\mathbb D^n, d\nu)$ consisting of analytic functions. In particular the space $A^2(\mathbb D^n)$ is a reproducing kernel Hilbert space, that is any $f\in A^2(\mathbb D^n)$ admits the representation
\Be
f(z)=P (f)(z)=\langle f, B (\cdot, z)\rangle=\int_{\mathbb D^n}f(w)B(z,w)d\nu(z),\,\,\,\textrm{for any}\,\,\,z\in \mathbb D^n,
\Ee
where the (weighted Bergman) kernel $B (\cdot,\cdot)$ is given by
\Be
B (z,w):=\frac{1}{\prod_{j=1}^n\left(1-\overline {w_j}z_j\right)^{2}}.
\Ee

$P$ is in fact the orthogonal projection from $L^2(\mathbb D^n)$ onto its closed subspace $A^2(\mathbb D^n)$ and it is called the Bergman projection. 

We denote by $\mathcal H(\mathbb D^n)$ the space of all analytic functions in $\mathbb D^n$. $H^\infty(\mathbb D^n)$ is the set of all bounded analytic functions in $\mathbb D^n$ that is $f\in H^\infty(\mathbb D^n)$ if $f\in \mathcal H(\mathbb D^n)$ and $$\|f\|_\infty:=\sup_{z\in \mathbb D^n}|f(z)|<\infty.$$

For $j\in \{1,\ldots, n\}$, we consider the operator $D_j$ defined for $f\in \mathcal{H}(\mathbb D^n)$ by
$$D_jf(z)=2f(z)+z_j\frac{\partial f}{\partial z_j}(z)=(2I+R_j)(f),$$ where $I$ stands for the identity operator. We put $D=D_1\ldots D_n$.
\vskip .2cm
The Bloch space of the polydisc $\mathbb D^n$ is denoted $\mathcal {B}(\mathbb D^n)$ and consists of all analytic functions $f$ such that
\Be\label{eq:Blochdefpoly}
\|f\|_{\mathcal {B}}:=|f(0)|+\sup_{z\in \mathbb D^n}\left[\prod_{j=1}^n(1-|z_j|^2)\right]\left| Df(z)\right|<\infty.
\Ee

The next space is called by several authors the Bloch space but here we name it the pointwise Bloch space to avoid any confusion with the Bloch space defined above.  A function $f$ analytic in $\mathbb D^n$ belongs to the pointwise Bloch space $\mathbb {B}(\mathbb D^n)$ if
\Be
\|f\|_{\mathbb {B}}:=|f(0)|+\sup_{z\in \mathbb D^n}\sum_{j=1}^n(1-|z_j|^2)\left|\frac{\partial f}{\partial z_j}(z)\right|<\infty.
\Ee
Note that a function is in the pointwise Bloch space if and only if it is a one parameter Bloch function in each variable. This justifies our choice of the name for this space.

Let us introduce some other spaces of analytic functions in $\mathbb D^n$. We start with the logarithmic Bloch space $\mathcal {B}_L(\mathbb D^n)$
which consists of analytic functions $f$ in $\mathbb D^n$ such that
\Be
\|f\|_{\mathcal {B}_L}:=|f(0)|+\sup_{z\in \mathbb D^n}\left[\prod_{j=1}^n(1-|z_j|^2)\log\frac{2}{1-|z_j|^2}\right]\left|Df(z)\right|<\infty.
\Ee

The above notion extends the notion of logarithmic Bloch space of the unit disc.

When it comes to the pointwise Bloch space, we also have a notion of logarithmic  Bloch space, $\mathbb B_L(\mathbb D^n)$. An analytic function $f$ belongs to $\mathbb B_L(\mathbb D^n)$ if
\Be
\|f\|_{\mathbb {B}_L}:=|f(0)|+\sup_{z\in \mathbb D^n}\left[\prod_{j=1}^n(1-|z_j|^2)\right]\left[\sum_{j=1}^n\log\frac{2}{1-|z_j|^2}\right]\left|Df(z)\right|<\infty.
\Ee

Our last space is the space ${\B}_{LL}(\mathbb D^n)$ defined as the space of all analytic functions $f$ in $\mathbb D^n$ such that there is a constant $C>0$ such that for any $K=\{k_1,\ldots, k_l\}\subseteq \{1,\cdots, n\}$,
 $$|f(0)|+\sup_{z\in \mathbb D^{n}}\left[\prod_{j\in K}(1-|z_j|^2)\log\frac{2}{1-|z_j|^2}\right]\left| D_{k_1}\ldots D_{k_l}f(z)\right|\le C.$$

The smallest constant in the above definition is denoted $\|f\|_{\mathcal B_{LL}}$.

Note that a function $f$ belongs to $\mathcal B_{LL}(\mathbb D^n)$ if and only if $f\in \mathcal B_L(\mathbb D^n)$ and for any $w\in \mathbb D^k$ fixed ($1\le k<n$), the function $f(\cdot,w)$ which is a function of $n-k$ variables, is uniformly in $\mathcal B_L(\mathbb D^{n-k})$.
\section{Statement of the results}\label{section3}
We recall that given $b\in A^2(\mathbb D^n)$, the (small) Hankel operator with symbol $b$,  $h_b$ is the operator defined for $f\in H^\infty(\mathbb D^n)$ by $$h_b(f):=P(b\overline {f}).$$
Our first result says that the space $\mathcal {B}_L(\mathbb D^n)$ is the exact range of symbols of bounded Hankel operators on $A^1(\mathbb D^n)$.
\btheo\label{theo:main01} Let $b\in \mathcal H(\mathbb D^n)$. Then the Hankel operator $h_b$ extends as a bounded operator on $A^1(\mathbb D^n)$ if and only if $b\in \mathcal B_L(\mathbb D^n)$.
\etheo

We also obtain that $\mathbb {B}_L(\mathbb D^n)$ is the exact range of symbols of bounded Hankel operators from $\mathbb {B}(\mathbb D^n)$ to $\mathcal {B}(\mathbb D^n)$.
\btheo\label{theo:main02} Let $b\in \mathcal H(\mathbb D^n)$. Then the Hankel operator $h_b$ extends as a bounded operator from $\mathbb B(\mathbb D^n)$ to $\mathcal B(\mathbb D^n)$ if and only if $b\in \mathbb B_L(\mathbb D^n)$.
\etheo

Given two Banach spaces of analytic functions $X$ and $Y$, the set of pointwise multipliers from $X$ to $Y$ is defined by
$$\mathcal {M}(X,Y):=\{g\in \mathcal H(\mathbb D^n):fg\in Y,\,\,\textrm{for any}\,\,\,f\in X\}.$$
When $X=Y$, we just write $\mathcal M(X)$ for $\mathcal M(X,X)$. The norm of the multiplication operator by $\phi$ from $X$ to $Y$ is denoted $\|M_\phi\|_{X\rightarrow Y}$ or $\|\phi\|_{X\rightarrow Y}$.

Our  first main result on pointwise multipliers is the following.
\btheo\label{theo:main1}
Let $\phi\in \mathcal H(\mathbb D^n)$. Then $\phi$ is a multiplier from $\mathbb B(\mathbb D^n)$ to $\mathcal B(\mathbb D^n)$ if and only if $\phi\in H^\infty(\mathbb D^n)\cap \mathbb B_L(\mathbb D^n)$. Moreover,
$$\|\phi\|_{\mathbb {B}\rightarrow \mathcal {B}}\backsimeq \|\phi\|_{\infty}+\|\phi\|_{\mathbb B_L}.$$
\etheo
Here is our characterization of the pointwise multipliers of $\mathcal B(\mathbb D^n)$.
\btheo\label{theo:main2}
Let $\phi\in \mathcal H(\mathbb D^n)$. Then $\phi$ is a multiplier of $\mathcal {B}(\mathbb D^n)$ if and only if
$$\phi\in H^\infty(\mathbb D^n)\cap \mathcal {B}_{LL}(\mathbb D^n).$$
Moreover,
$$\|\phi\|_{\mathcal {B}\rightarrow \mathcal {B}}\backsimeq \|\phi\|_{\infty}+\|\phi\|_{\mathcal B_{LL}}.$$
\etheo

 In Section \ref{section4}, we give some useful properties of $\mathcal B(\mathbb D^n)$ and their logarithmic counterparts.  The proofs of Theorem \ref{theo:main01} and Theorem \ref{theo:main02} are given in Section \ref{section5}.  In  Section \ref{section6} we prove our results on multipliers from $\mathbb B(\mathbb D^n)$ to $\mathcal B(\mathbb D^n)$ and the multipliers of $\mathcal B(\mathbb D^n)$.  In the last section of this paper, we add some comments and remarks.

Finally, all over the text, $C$ will be a constant not necessarily the same at each occurrence. We will also use the notation $C_k$
 to express the fact that the constant depends on the underlined parameter. Given two positive quantities $A$ and $B$, the notation
 $A\lesssim B$ means that $A\le CB$ for some positive constant $C$. When $A\lesssim B$ and $B\lesssim A$, we write $A\backsimeq B$.

\section{Useful results on the Bloch spaces}\label{section4}

\subsection{The Bloch space of the unit disc}
On the unit disc $\mathbb D$ of the complex plane $\mathbb C$, the Bloch space $\mathcal B=\mathcal {B}(\mathbb D)$ consists of analytic functions $f$ such that
\Be\label{eq:Blochdefdisc}
\sup_{z\in \mathbb D}(1-|z|^2)|f'(z)|<\infty.
\Ee
The following norm makes $\mathcal B(\mathbb D)$ a Banach space:
\Be\label{eq:Blochdefdiscnorm}
\|f\|_{\mathcal B}:=|f(0)|+\sup_{z\in \mathbb D}(1-|z|^2)|f'(z)|<\infty.
\Ee

The Bloch space strictly contains the space $H^\infty(\mathbb D)$. A typical example of function in the Bloch space that does not belong to $H^\infty(\mathbb D)$, is the function $f(z)=\log (1-\overline {a}z)$ which is uniformly in $\mathcal B(\mathbb D)$. That is its $\mathcal B$-norm is bounded by a constant that does not depend on the complex number $a$ (see \cite{benoit1}).

Note that in the above definition, $f'(z)$ can be replaced by $Rf(z)=zf'(z)$ and that equivalent norms are obtained by considering any derivative of higher order. That is for any integer $k\ge 1$,
\Beas
\|f\|_{\mathcal B} &\simeq& |f(0)|+\sup_{z\in \mathbb D}(1-|z|^2)^k|f^{(k)}(z)|\\ &\simeq& |f(0)|+\sup_{z\in \mathbb D}(1-|z|^2)^k|R^kf(z)|.
\Eeas

In general, higher order derivatives can be replaced by the so-called fractional derivatives (for a definition, see \cite{Zhu}). In particular, if for the integer $k\ge 1$ we define the operator $$D^k=[(1+k)I+R]\cdots [2I+R]$$ where $I$ is the identity, then
\Be
\|f\|_{\mathcal B} \simeq |f(0)|+\sup_{z\in \mathbb D}(1-|z|^2)^k|D^kf(z)|
\Ee
(see \cite{BL}).
Let us still denote by $P$ the orthogonal projection from $L^2(\mathbb D)$ onto its closed subspace $A^2(\mathbb D)$. For $b\in A^2(\mathbb D)$, we densely
defined the (small) Hankel operator with symbol $b$ on $A^2(\mathbb D)$ by $$h_b(f):=P(b\overline {f}).$$
There are some other equivalent characterizations of the Bloch space (see \cite{BL, Zhu}).
\bprop Let $b$ be an analytic function in the unit discc $\mathbb D$. Then the following are equivalent.
\begin{itemize}
\item[(i)] $b\in \mathcal B(\mathbb D)$;
\item[(ii)] $b=Pg$ for some $g\in L^\infty(\mathbb D)$;
\item[(iii)] $b$ belongs to the dual space $\left(A^1(\mathbb D)\right)^*$ of $A^1(\mathbb D)$ under the pairing
$$\langle f,g\rangle:=\lim_{r\rightarrow 1}\int_{\mathbb D}f(rz)\overline {g(z)}d\nu(z);$$
\item[(iv)] The Hankel operator $h_b$ is bounded on $A^2(\mathbb D)$.
\end{itemize}
\eprop

To deal with multipliers of $\mathcal B(\mathbb D)$, we recall that the logarithmic counterpart of the Bloch space called the logarithmic Bloch space and denoted $\mathcal B_L(\mathbb D)$, consists of all holomorphic functions $f$ in $\mathbb D$ such that
\Be
\|f\|_{\mathcal B_L} \simeq |f(0)|+\sup_{z\in \mathbb D}(1-|z|^2)\left(\log\frac{2}{1-|z|^2}\right)|f'(z)|<\infty.
\Ee
Remark that the same observations made on $\mathcal B$ about higher order derivatives also work for the logarithmic Bloch space (see for example \cite{BL}). That is
\Beas
\|f\|_{\mathcal B_L} &\simeq& |f(0)|+\sup_{z\in \mathbb D}(1-|z|^2)^k\left(\log\frac{2}{1-|z|^2}\right)|f^{(k)}(z)|\\ &\simeq& |f(0)|+\sup_{z\in \mathbb D}(1-|z|^2)^k\left(\log\frac{2}{1-|z|^2}\right)|R^kf(z)|\\ &\simeq& |f(0)|+\sup_{z\in \mathbb D}(1-|z|^2)^k\left(\log\frac{2}{1-|z|^2}\right)|D^kf(z)|.
\Eeas

One has the following characterization of the multiplier algebra of the Bloch space in the unit disc.
\bprop\label{prop:multionepara}
Let $b\in \mathcal {H}(\mathbb D)$. Then $b\in \mathcal {M}(\mathcal B(\mathbb D))$ if and only if $$b\in H^\infty(\mathbb D)\cap \mathcal B_L(\mathbb D).$$ Moreover,
$$\|M_b\|_{\mathcal {B}\rightarrow \mathcal B}\simeq \|b\|_{\infty}+\|b\|_{\mathcal B_L}.$$
\eprop
\begin{proof}
The proof is quite standard, we give it here as it will guide us along the text. One thing that one needs to know to prove the sufficiency in the above proposition, is the following pointwise estimate of functions in the Bloch space.
$$|f(z)|\le C\left(\log\frac{2}{1-|z|^2}\right)\|f\|_{\mathcal B}.$$
The above pointwise estimate can be combined with the fact that the function $\log(1-\overline {a}z)$ is uniformly in $\mathcal B(\mathbb D)$ to prove that any element of $\mathcal {M}(\mathcal B(\mathbb D))$ is bounded (see \cite{benoit1, zhao}).

\vskip .2cm
Suppose that $b$ satisfies the condition in the proposition. Then for any $f\in \mathcal B(\mathbb D)$ and any $z\in \mathbb D$,
\Beas
(1-|z|^2)|\left(fb\right)'(z)| &=& (1-|z|^2)|f(z)b'(z)+f'(z)b(z)|\\ &\le& (1-|z|^2)|f(z)||b'(z)|+(1-|z|^2)|f'(z)||b(z)|\\ &\le&  C(1-|z|^2)||b'(z)|\left(\log\frac{2}{1-|z|^2}\right)\|f\|_{\mathcal B}+\|b\|_{\infty}(1-|z|^2)|f'(z)|\\ &\le& C\|f\|_{\mathcal B}\left(\|b\|_{\infty}+\|b\|_{\mathcal B_L}\right).
\Eeas

Now suppose that $b$ is multiplier of $\mathcal B(\mathbb D)$. That is there exists a constant $C>0$ such that for any $f\in \mathcal B(\mathbb D)$ and $z\in \mathbb D$,
$$(1-|z|^2)|\left(fb\right)'(z)| = (1-|z|^2)|f(z)b'(z)+f'(z)b(z)|\le C\|f\|_{\mathcal B}.$$
We test this inequality with the function $f(z)=f_a(z)=\log(1-\overline {a}z)$, $a\in \mathbb D$ fixed. It comes that
$$(1-|z|^2)|\log(1-\overline {a}z)b'(z)+\frac{\overline {a}}{1-\overline {a}z}b(z)|\le C.$$
Putting $z=a$, it comes since $b\in H^\infty(\mathbb D)$ that for any $z\in \mathbb D$,
$$(1-|z|^2)\left(\log\frac{2}{1-|z|^2}\right)|b'(z)|\le C<\infty.$$
The proof is complete.
\end{proof}
\subsection{The product Bloch space and its logarithmic counterpart}\label{section5}

For $K=(k_1,\ldots,k_n)$, $k_j\in \mathbb N$, we define on $\mathcal H(\mathbb D^n)$ the operator $D^K$ given by
$$D^K=D_1^{k_1}\ldots D_n^{k_n},$$
where $D_j^{k_j}f(z)=[(1+k_j)I+R_j]\ldots [2I+R_j]$ with $R_jf(z)=z_j\frac{\partial f}{\partial z_j}(z)$. 

As in the one parameter setting, we obtain equivalent norms on $\mathcal {B}(\mathbb D^n)$ by using higher order derivatives in each variable. More precisely, for a vector $K=(k_1,\cdots,k_n)\in \mathbb N^n$ , we have
\Be\label{eq:equivBloch}\|f\|_{\mathcal B}\simeq |f(0)|+\sup_{z\in \mathbb D^n}\left[\prod_{j=1}^n(1-|z_j|^2)^{k_j}\right]|D^Kf(z)|\Ee
Also,
$$\|f\|_{\mathcal B_L}\backsimeq |f(0)|+\sup_{z\in \mathbb D^n}\left[\prod_{j=1}^n(1-|z_j|^2)^{k_j}\log\frac{2}{1-|z_j|^2}\right]|D^Kf(z)|$$
(see also \cite{Har, HarLus}).
\vskip .2cm

Let us observe the following formula that can be proved as in the one parameter situation.
\blem\label{lem:integrationbyparts}
Let $f$ and $g$ be two analytic polynomials in $\mathbb D^n$. Then for any $K=(k_1,\cdots,k_n)\in \mathbb N_0^n$, there exists a constant $C=C_{K,n}$ such that the following formula holds
\Be\label{eq:integrationbyparts}
\int_{\mathbb D^n}f(z)\overline {g(z)}d\nu(z)=C\int_{\mathbb D^n}f(z)(1-||z||^2)^K\overline {D^Kg(z)}d\nu(z),
\Ee
$$(1-||z||^2)^K:=\prod_{j=1}^n(1-|z_j|^2)^{k_j}.$$
\elem

The following first equivalent characterization of the Bloch space was obtained by K. Zhu in \cite{Zhu1}.
\blem\label{lem:equivblochprod1}
Let $f$ be holomorphic in $\mathbb D^n$. Then the following assertions are equivalent.
\begin{itemize}
\item[(i)] $f\in \mathcal {B}(\mathbb D^n)$.
\item[(ii)] There exists a function $g\in L^\infty(\mathbb D^n)$ such that
$$f(z)=\int_{\mathbb D^n}\frac{g(w)}{\prod_{j=1}^n\left(1-z_j\overline {w_j}\right)^2}d\nu(z),\,\,\,z\in \mathbb D^n.$$
Moreover, $\|f\|_{\mathcal B}\backsimeq \|g\|_\infty$.
\end{itemize}
\elem

We refer also to \cite{Constantin} for the following duality result which provides another equivalent definition of $\mathcal B(\mathbb D^n)$.
\blem\label{lem:duality}
The topological dual space $\left(A^1(\mathbb D^n)\right)^*$ of $A^1(\mathbb D^n)$ identifies with $\mathcal B(\mathbb D^n)$ under the duality pairing
\Be\label{eq:dualitypairing}
\langle f,g\rangle:=\lim_{r\rightarrow {\bf 1}}\int_{\mathbb D^n}f(rz)\overline {g(z)}d\nu (z), f\in A^1(\mathbb D^n),g\in \mathcal B(\mathbb D^n).
\Ee
\elem

The following pointwise estimate of functions in $\mathcal B(\mathbb D^n)$ will be useful, it can be obtained from the representation of Bloch functions in Lemma \ref{lem:equivblochprod1}.
\blem There exists a constant $C>0$ such that for any $f\in \mathcal B(\mathbb D^n)$ and any $z\in \mathbb D^n$,
\Be\label{eq:pointprodbloch}
|f(z)|\le C\left(\prod_{j=1}^n\log\frac{2}{1-|z_j|^2}\right)\|f\|_{\mathcal B}.
\Ee
\elem

Another equivalent definition of the Bloch space is in term of symbols of bounded Hankel operators on the Bergman space $A^p(\mathbb D^n)$, $1<p<\infty$.
\bprop Let $1<p<\infty$ and $b\in \mathcal {H}(\mathbb D^n)$. Then the Hankel operator $h_b$ extends as a bounded operator on $A^p(\mathbb D^n)$ if and only if $b=Pg$ for some $g\in L^\infty(\mathbb D^n)$.
\eprop
\begin{proof}
The proof of this result can be found in \cite{Constantin}. We give an alternative proof here using test functions. Let us first suppose that $b$ is as above, then $b\in \mathcal B(\mathbb D^n)$ by Lemma \ref{lem:equivblochprod1}. Hence, for any $f\in A^p(\mathbb D^n)$ and $m\in A^q(\mathbb D^n)$, $pq=p+q$, using the duality in Lemma \ref{lem:duality}, we obtain
$$|\langle h_b(f),m\rangle|=|\langle P(b\bar {f}),m\rangle|=|\langle b,fm\rangle|\le \|b\|_{\mathcal B(\mathbb D^n)}\|fm\|_1\le \|b\|_{\mathcal B(\mathbb D^n)}\|f\|_p\|m\|_q.$$
That is $h_b$ is bounded on $A^p(\mathbb D^n)$ as soon as $b$ is as in the statement of the proposition.

Now suppose that $b$ is analytic on $\mathbb D^n$ and such that $h_b$ extends as a bounded operator on $A^p(\mathbb D^n)$. Note that this means that there is a constant $C>0$ such that for any $f\in A^p(\mathbb D^n)$ and any $g\in A^q(\mathbb D^n)$, $pq=p+q$,
\Be\label{eq:necesshank}
|\langle b,fg\rangle|\le C\|f\|_p\|g\|_q
\Ee
Let $a\in \mathbb D^n$ be fixed and put $$f(z)=f_a(z)=\prod_{j=1}^n\frac{(1-|a_j|^2)^{k_j/p}}{(1-z_j\bar {a_j})^{(2+k_j)/p}}$$ and
$$g(z)=g_a(z)=\prod_{j=1}^n\frac{(1-|a_j|^2)^{k_j/q}}{(1-z_j\bar {a_j})^{(2+k_j)/q}},\,\,\, k_j\in \N, \,\,\,j=1,2,\cdots,n.$$ Observe that $f$ and $g$ are uniformly in $A^p(\mathbb D^n)$ and $A^q(\mathbb D^n)$ respectively. Taking $f=f_a$ and $g=g_a$ in (\ref{eq:necesshank}), we obtain
\Beas
 C\ge |\langle b,fg\rangle|=\left(\prod_{j=1}^n(1-|a_j|^2)^{k_j}\right)\left|\lim_{r\rightarrow {\bf 1}}\int_{\mathbb D^n}\frac{b(rw)}{\prod_{j=1}^n(1-a_j\bar {w_j})^{2+k_j}}d\nu(w)\right|.
\Eeas
That is for any $a\in \mathbb D^n$ and any $K=(k_1,\ldots,k_n)\in \N^n$, $$\left(\prod_{j=1}^n(1-|a_j|^2)^{k_j}\right)|D^Kb(a)|\le C$$
and consequently, $b\in \mathcal B(\mathbb D^n)$. The conclusion then follows from Lemma \ref{lem:equivblochprod1}.
\end{proof}
\begin{remark} As observed in \cite{Constantin}, the above result is equivalent to saying that any $f\in A^1(\mathbb D^n)$ admits a representation of the form
$$f(z)=\sum_{j}f_j(z)g_j(z),\,\,z\in \mathbb D^n, f_j\in A^p(\mathbb{D}), g_j\in A^q(\mathbb{D}),$$ with $$\sum_{j\in \N_0}\|f_j\|_p\|g_j\|_q\le C\|f\|_1,\,\,pq=p+q,\,1<p<\infty.$$
\end{remark}
\subsection{The pointwise Bloch space}\label{section6}
Let us start by considering some observations.
\blem\label{lem:pointwisesmallbloch}
The following assertions hold
\begin{itemize}
\item[(1)]
The function $f(z)=\sum_{j=1}^nf_j(z_j)$, with $f_j\in \mathcal B(\mathbb D)$ belongs  to $\mathbb B(\mathbb D^n)$ and $$\|f\|_{\mathcal B(\mathbb D^n)}\le \sum_{j=1}^n\|f_j\|_{\mathcal B(\mathbb D)}.$$
\item[(2)] There is a constant $C>0$ such that for any $f\in \mathbb B(\mathbb D^n)$ and any $z\in \mathbb D^n$,
\Be\label{eq:pointwisesmallbloch}
|f(z)|\le C\left[\sum_{j=1}^n\log\frac{2}{1-|z_j|^2}\right]\|f\|_{\mathbb B}
\Ee
and this is sharp.
\end{itemize}
\elem
\begin{proof}
Assertion $(1)$ is direct from the definition of $\mathbb B(\mathbb D^n)$. Let us prove $(2)$. We have for any $z=(z_1,\ldots, z_n)\in \mathbb D^n$,
\Beas
f(z)-f(0) &=& \int_0^1\frac{d f(sz)}{ds}ds\\ &=& \sum_{j=1}^n\int_0^1 z_j\frac{\partial f}{\partial z_j}(sz)ds.
\Eeas
It easily follows using the definition of $\mathbb B(\mathbb D^n)$ that
\Beas
|f(z)| &\le& |f(0)|+\sum_{j=1}^n\int_0^1 |z_j|\left|\frac{\partial f}{\partial z_j}(sz)\right|ds\\
&\le& |f(0)|+\|f\|_{\mathcal B(\mathbb D^n)}\sum_{j=1}^n\int_0^1\frac{|z_j|}{1-s^2|z_j|^2}ds\\ &\le& C\left(\sum_{j=1}^n\log\frac{2}{1-|z_j|^2}\right)\|f\|_{\mathcal B}.
\Eeas

Sharpness follows by testing (\ref{eq:pointwisesmallbloch}) with the function $f(z)=f_a(z)=\sum_{j=1}^n\log(1-\overline {a_j}z_j)$.
\end{proof}

\begin{remark}: Let us observe that $\mathbb B(\mathbb D^n)$ is a strict subspace of $\mathcal B(\mathbb D^n)$. To see this, one only needs to observe that for $f_j\in \mathcal {B}(\mathbb D)$, $j=1,\cdots,n$, the tensor product $$\left(f_1\otimes f_2\otimes\cdots\otimes f_n\right)(z_1,\cdots,z_n)=\prod_{j=1}^nf_j(z_j)$$ belongs to  $\mathcal {B}(\mathbb D^n)$ with
$$\|f_1\otimes f_2\otimes\cdots\otimes f_n\|_{\mathcal {B}(\mathbb D^n)}\le \prod_{j=1}^n\|f_j\|_{\mathcal {B}(\mathbb D)}$$
while $f_1\otimes f_2\otimes\cdots\otimes f_n$ belongs to $\mathbb B(\mathbb D^n)$ only if each $f_j$ belongs to $H^\infty(\mathbb D)$, $j=1,2,\cdots,n$.
\end{remark}
\vskip .2cm
 To $\vec {j}=({j_1},\ldots,{j_l})$, $1\le l\le n$, we associate the set $J=\{j_1,\ldots,j_l\}\subseteq \{1,\ldots,n\}$. We denote by $D_{\vec j}$ the differential operator defined by
$$D_{\vec j}f(z)=D_{j_1}\ldots D_{j_l}f(z).$$
We observe the following.
\blem\label{lem:smallblochpartialBloch}
Let $\vec {j}=({j_1},\ldots,{j_l})$, $1\le l\le n$, $J=\{j_1,\ldots,j_l\}\subseteq \{1,\ldots,n\}$ be given. Then for any $f\in \mathbb B(\mathbb D^n)$,
\Be\label{eq:smallblochpartialBloch}\sup_{z\in \mathbb D^n}\left[\prod_{k\in J}(1-|z_k|^2)\right]|D_{\vec j}f(z)|\lesssim \|f\|_{\mathbb B}.\Ee
\elem
\begin{proof}
We can suppose without loss of generality that $\vec j=(1,\ldots,l)$, $1\le l\le n$, so that $J=\{1,\ldots, l\}$. For $z=(z_1,\ldots,z_n)$, we set $w=(z_{l+1},\ldots,z_n)\in \mathbb D^{n-l}$. We observe that for $w\in \mathbb D^{n-l}$ fixed, the function $g=f(\cdot,w)$ is uniformly in $\mathbb B(\mathbb D^l)$ whenever $f\in \mathbb B(\mathbb D^n)$ with $\|g\|_{\mathbb B(\mathbb D^l)}\le \|f\|_{\mathbb B(\mathbb D^n)}$. Hence, $g=f(\cdot, w)$ is uniformly in $\mathcal B(\mathbb D^l)$ with
$$\sup_{a\in \mathbb D^l}\left[\prod_{k=1}^l(1-|a_k|^2)\right]|D_{\vec j}f(a,w)|\le \|g\|_{\mathbb B(\mathbb D^l)}\le \|f\|_{\mathbb B(\mathbb D^n)}.$$
This proves that (\ref{eq:smallblochpartialBloch}) holds.
\end{proof}
\vskip .2cm
Let us observe that as for the space $\mathcal B_L(\mathbb D^n)$, we have the following equivalent definition for the  logarithmic Bloch-type space that we have denoted $\mathbb B_L(\mathbb D^n)$. 
\Beas \|f\|_{\mathbb B_L} \backsimeq |f(0)|+\sup_{z\in \mathbb D^n}\left[\prod_{j=1}^n(1-|z_j|^2)^{k_j}\right]\left[\sum_{j=1}^n\log\frac{2}{1-|z_j|^2}\right]|D^Kf(z)|<\infty,
\Eeas
$K=(k_1,\cdots,k_n)\in \mathbb {N}^n$.

\section{Endpoint results for Hankel operators}\label{section7}
We now prove Theorem \ref{theo:main01}. This provides  an equivalent definition of $\mathcal B_L(\mathbb D^n)$ in terms of symbols of bounded Hankel operators on $A^1(\mathbb D^n)$ as in the one parameter case (see for example \cite{BL}). We will then be calling $\mathcal B_L(\mathbb D^n)$ the product logarithmic Bloch space.
\begin{proof}[Proof of Theorem \ref{theo:main01}]
We start by the easy part which is the sufficiency. What we would like to prove is that giving any $f\in A^1(\mathbb D^n)$, the function $h_b(f)$ belongs to $A^1(\mathbb D^n)$ under the condition that $b\in \mathcal B_L(\mathbb D^n)$ or equivalently that
$$\|h_b(f)\|_1\le C\|f\|_1\,\,\,\textrm{with}\,\,\, C \,\,\,\textrm{not depending on}\,\,\,f.$$
To prove this, we observe with the help of Lemma \ref{lem:integrationbyparts} that we have in particular that
\Beas h_b(f)(z) &=& \int_{\mathbb D^n}\frac{b(w)\overline {f(w)}}{\prod_{j=1}^n(1-\overline {w_j}z_j)^2}d\nu(w)\\ &=& C_{n,K}\int_{\mathbb D^n}\frac{[(1-\|w\|^2)^KD^Kb(w)]\overline {f(w)}}{\prod_{j=1}^n(1-\overline {w_j}z_j)^2}d\nu(w),
\Eeas
for any $K\in \mathbb N^n$. It follows using this observation and \cite[Proposition 1.4.10]{Rudin} that
\Beas
\|h_b(f)\|_1 &=& \int_{\mathbb D^n}\left|h_b(f)(z)\right|d\nu(z)\\ &\le& C\int_{\mathbb D^n}\left|\int_{\mathbb D^n}\frac{[(1-\|w\|^2)^KD^Kb(w)]\overline {f(w)}}{\prod_{j=1}^n(1-\overline {w_j}z_j)^2}d\nu(w)\right|d\nu(z)\\ &\le& C\int_{\mathbb D^n}(1-\|w\|^2)^K|D^Kb(w)||f(w)|\left(\int_{\mathbb D^n}\frac{d\nu(z)}{\prod_{j=1}^n|1-\overline {w_j}z_j|^2}\right)d\nu(w)\\ &\le& C\int_{\mathbb D^n}(1-\|w\|^2)^K\left(\prod_{j=1}^n\log\frac{2}{1-|w_j|^2}\right)|D^Kb(w)||f(w)|d\nu(w)\\ &\le& C\|b\|_{\mathcal B_L}\int_{\mathbb D^n}|f(w)|d\nu(w)\\ &=& C\|b\|_{\mathcal B_L}\|f\|_1.
\Eeas

The converse is equivalent to saying that if $b$ is such that there exists a constant $C>0$ so that for any $f\in A^1(\mathbb D^n)$ and any $g\in \mathcal B(\mathbb D^n)$,
\Be\label{eq:testhank}
\left|\langle b,fg\rangle\right|\le C\|f\|_1\|g\|_{\mathcal B},
\Ee
then $b\in \mathcal B_L(\mathbb D^n)$.

We will need the following lemma.
\blem\label{lem:logprodlemma} Let $w_j, z_j\in \mathbb D$, $j=1,\ldots, n$ be given. The following formula holds.
\Be\label{eq:logfacto}
\prod_{j=1}^n\log(1-\overline {w_j}z_j)=\sum_{L\subseteq \{1,\ldots, n\}}(-1)^{[L]+1}T_L+\prod_{j=1}^n\log(1-|w_j|^2)
\Ee
where $$T_L=\left(\prod_{l\in L}\log\frac{1-\overline {w_l}z_l}{1-|w_l|^2}\right)\left(\prod_{l\in {^cL}}\log(1-\overline {w_l}z_l)\right),$$
$[L]$ being the cardinality of the set $L$ and $^cL$ its complementary in $\{1,\ldots, n\}$.
\elem
\begin{proof}
One easily checks that for $n=2$, we have
\Beas \prod_{j=1}^2\log(1-\overline {w_j}z_j)&=&  \log\frac{1-\overline {w_1}z_1}{1-|w_1|^2}\log(1-\overline {w_2}z_2)+\log\frac{1-\overline {w_2}z_2}{1-|w_2|^2}\log(1-\overline {w_1}z_1)\\ & -& \log\frac{1-\overline {w_1}z_1}{1-|w_1|^2}\log\frac{1-\overline {w_2}z_2}{1-|w_2|^2}+\log(1-|w_1|^2)\log(1-|w_2|^2).
\Eeas
Next we suppose that (\ref{eq:logfacto}) holds for $n\ge 2$ and prove that it then also holds for $n+1$.
Using our hypothesis, we obtain
\Beas
\prod_{j=1}^{n+1}\log(1-\overline {w_j}z_j) &=& \left(\prod_{j=1}^n\log(1-\overline {w_j}z_j)\right)\log(1-\overline {w_{n+1}}z_{n+1})\\ &=& T_1+\prod_{j=1}^n\log(1-|w_j|^2)\log(1-\overline {w_{n+1}}z_{n+1}).
\Eeas
where $$T_1=\sum_{L\subseteq \{1,\cdots, n\}}(-1)^{[L]+1}\left(\prod_{l\in L}\log\frac{1-\overline {w_l}z_l}{1-|w_l|^2}\right)\left(\prod_{l\in {^cL}}\log(1-\overline {w_l}z_l)\right)\log(1-\overline {w_{n+1}}z_{n+1})$$
Before going ahead, let us deal with the second term in the sum on the right hand side of the above equality. We clearly have using our hypothesis again that
\Beas
M &:=& \prod_{j=1}^n\log(1-|w_j|^2)\log(1-\overline {w_{n+1}}z_{n+1})\\ &=& \prod_{j=1}^n\log(1-|w_j|^2)\log(\frac{1-\overline {w_{n+1}}z_{n+1}}{1-|w_{n+1}|^2})+\prod_{j=1}^{n+1}\log(1-|w_j|^2)\\ &=& T_2 + T_3 + \prod_{j=1}^{n+1}\log(1-|w_j|^2),
\Eeas

where
$$T_2=\left(\prod_{j=1}^n\log(1-\overline {w_j}z_j)\right)\log(\frac{1-\overline {w_{n+1}}z_{n+1}}{1-|w_{n+1}|^2})$$
and 
$$T_3=(-1)^{n}\left(\prod_{l=1}^n\log\frac{1-\overline {w_l}z_l}{1-|w_l|^2}\right)\log(\frac{1-\overline {w_{n+1}}z_{n+1}}{1-|w_{n+1}|^2})$$
Taking this into the expansion of $\prod_{j=1}^{n+1}\log(1-\overline {w_j}z_j)$, we obtain
\Beas
\prod_{j=1}^{n+1}\log(1-\overline {w_j}z_j) &=& T_1+T_2+ T_3 +\prod_{j=1}^{n+1}\log(1-|w_j|^2)\\ &=& \sum_{L\subseteq \{1,\ldots, n+1\}}(-1)^{[L]+1}\left(\prod_{l\in L}\log\frac{1-\overline {w_l}z_l}{1-|w_l|^2}\right)\left(\prod_{l\in {^cL}}\log(1-\overline {w_l}z_l)\right)\\ &+& \prod_{j=1}^{n+1}\log(1-|w_j|^2).
\Eeas
The proof of the lemma is complete.
\end{proof}
Coming back to the proof of the necessity part of the theorem, we test (\ref{eq:testhank}) with $$f(z)=f_w(z)=\prod_{j=1}^n\frac{(1-|w_j|^2)^{k_j}}{\left(1-\overline {w_j}z_j\right)^{2+k_j}}\,\,\,k_j\in \mathbb N_0,$$
and $$g(z)=g_w(z)=\prod_{j=1}^{n}\log(1-\overline {w_j}z_j).$$
Clearly $f$ and $g$ are uniformly in $A^1(\mathbb D^n)$ and $\mathcal B(\mathbb D^n)$ respectively.

Next, we take $f_w$ and the expansion of $g_w$ obtained in Lemma \ref{lem:logprodlemma} into (\ref{eq:testhank}) to obtain that
$$C\ge |\langle b, fg\rangle|=\left|\left(\prod_{j=1}^{n}\log(1-|w_j|^2)\right)\lim_{r\rightarrow {\bf 1}}\int_{\mathbb D^n}\prod_{j=1}^n\frac{(1-|w_j|^2)^{k_j}}{\left(1-\overline {w_j}z_j\right)^{2+k_j}}\overline {b(rz)}d\nu(z)+T\right|$$
where writing $$G_w^L(z)=\left(\prod_{l\in L}\log\frac{1-\overline {w_l}z_l}{1-|w_l|^2}\right) \left(\prod_{l\in {^cL}}\log(1-\overline {w_l}z_l)\right),$$
\Beas
 T &=& \lim_{r\rightarrow {\bf 1}}\int_{\mathbb D^n}\left(\sum_{L\subseteq \{1,\cdots, n\}}(-1)^{[L]+1}G_w^L(z)\right)\left(\prod_{j=1}^n\frac{(1-|w_j|^2)^{k_j}}{\left(1-\overline {w_j}z_j\right)^{2+k_j}}\right)\overline {b(rz)}d\nu(z)\\ &=& \sum_{L\subseteq \{1,\cdots, n\}}(-1)^{[L]+1}\lim_{r\rightarrow {\bf 1}}\int_{\mathbb D^n}f_w^L(z)g_w^L(z)\overline {b(rz)}d\nu(z)
\Eeas
where
\Beas f_w^L(z) &:=& \left(\prod_{j=1}^n\frac{(1-|w_j|^2)^{k_j}}{\left(1-\overline {w_j}z_j\right)^{2+k_j}}\right)\left(\prod_{l\in L}\log\frac{1-\overline {w_l}z_l}{1-|w_l|^2}\right)\\ &=& \left(\prod_{j\in {^cL}}\frac{(1-|w_j|^2)^{k_j}}{\left(1-\overline {w_j}z_j\right)^{2+k_j}}\right)\left(\prod_{j\in L}\frac{(1-|w_j|^2)^{k_j}}{\left(1-\overline {w_j}z_j\right)^{2+k_j}}\log\frac{1-\overline {w_j}z_j}{1-|w_j|^2}\right)
\Eeas
and
$$g_w^L(z):=\prod_{l\in {^cL}}\log(1-\overline {w_l}z_l).$$
Clearly, $g_w^L$ is uniformly in $\mathcal B(\mathbb D^n)$. Observing that in the unit disc of $\mathbb C$, the function
$\frac{(1-|w_j|^2)^{k_j}}{\left(1-\overline {w_j}z_j\right)^{2+k_j}}\log\frac{1-\overline {w_j}z_j}{1-|w_j|^2}$ is uniformly in $A^1(\mathbb D)$ (see \cite{BL}), we conclude that $f_w^L$ is also uniformly in $A^1(\mathbb D^n)$. Hence, applying (\ref{eq:testhank}) to $f_w^L$ and $g_w^L$, we obtain that
$$|T|\le \sum_{L\subseteq \{1,\cdots, n\}}\left|\langle b, f_w^Lg_w^L\rangle\right| \le C.$$
We deduce that
$$\left(\prod_{j=1}^{n}(1-|w_j|^2)^{k_j}\left|\log(1-|w_j|^2)\right|\right)\left|\lim_{r\rightarrow {\bf 1}}\int_{\mathbb D^n}\frac{\overline {b(rz)}}{\prod_{j=1}^n\left(1-\overline {w_j}z_j\right)^{2+k_j}}d\nu(z)\right|\le C.$$
That is $$\left(\prod_{j=1}^{n}(1-|w_j|^2)^{k_j}\left|\log(1-|w_j|^2)\right|\right)|D^Kb(w)|\le C<\infty, \,\,\,\textrm{for any}\,\,\,w\in \mathbb D^n.$$
We conclude that $b\in \mathcal B_L(\mathbb D^n)$. The proof is complete.
\end{proof}

Let us now prove Theorem \ref{theo:main02} that provides an equivalent definition of $\mathbb B_L(\mathbb D)$ in terms of symbols of bounded Hankel operators.
\begin{proof}[Proof of Theorem \ref{theo:main02}]
Let us start by the sufficiency. Assume $b\in \mathbb B_L(\mathbb D^n)$. Then for any $f\in \mathbb {B}(\mathbb D^n)$ and any $g\in A^1(\mathbb D^n)$, and for $K=(k_1,\cdots,k_n)\in \N^n$,
\Beas
|\langle h_b(f),g\rangle| &=& |\langle b,fg\rangle|\\  &\le& C\lim_{r\rightarrow {\bf 1}}\int_{\mathbb D^n}\left|[(1-\|w\|^2)^KD^Kb(w)]\overline {f(rw)g(rw)}\right|d\nu(z)\\ &\le& C\lim_{r\rightarrow {\bf 1}}\int_{\mathbb D^n}(1-\|w\|^2)^K|D^Kb(w)||f(rw)||g(rw)d\nu(w)\\ &\le& C\|f\|_{\mathbb B(\mathbb D^n)}\lim_{r\rightarrow {\bf 1}}\int_{\mathbb D^n}(1-\|w\|^2)^K\left(\sum_{j=1}^n\log\frac{2}{1-|w_j|^2}\right)|D^Kb(w)||g(rw)|d\nu(w)\\ &\le& C\|f\|_{\mathbb B(\mathbb D^n)}\|b\|_{\mathbb B_L(\mathbb D^n)}\int_{\mathbb D^n}|g(w)|d\nu(w)\\ &=& C\|f\|_{\mathbb B(\mathbb D^n)}\|b\|_{\mathbb B_L(\mathbb D^n)}\|g\|_1.
\Eeas
Thus $$\|h_b(f)\|_{\mathcal B(\mathbb D^n)}=\sup_{g\in A^1(\mathbb D^n), \|g\|_1\le 1}|\langle h_b(f),g\rangle|\le C\|f\|_{\mathbb B(\mathbb D^n)}\|b\|_{\mathbb B_L(\mathbb D^n)}.$$
That is $h_b$ is bounded from $\mathbb {B}(\mathbb D^n)$ to $\mathcal {B}(\mathbb D^n)$ for any $b\in \mathbb B_L(\mathbb D^n)$.
\vskip .2cm
For the converse, we have to prove that if $b$ is such that there exists a constant $C>0$ so that for any $f\in \mathbb B(\mathbb D^n)$, and any $g\in A^1(\mathbb D^n)$,
\Be\label{eq:testhank1}
\left|\langle b,fg\rangle\right|\le C\|f\|_{\mathbb {B}(\mathbb D^n)}\|g\|_1,
\Ee
then $b\in \mathbb B_L(\mathbb D^n)$. For this, we test (\ref{eq:testhank1}) with
$$f_a(z)=\sum_{j=1}^n\log(1-\overline {a_j}z_j),\,\,\,a=(a_1,a_2,\cdots,a_n)\in \mathbb D^n,$$
and $$g_a(z)=\prod_{j=1}^n\frac{(1-|a_j|^2)^{k_j}}{\left(1-\overline {a_j}z_j\right)^{2+k_j}}\,\,\,k_j\in \mathbb N_0.$$
We obtain at the first step using that both $f_a$ and $g_a$ are uniformly in $\mathbb B(\mathbb D^n)$ and $A^1(\mathbb D^n)$ respectively,
\Be\label{eq:step1}
C\ge \left|\langle b,fg\rangle\right| = \left|\sum_{l=1}^n\lim_{r\rightarrow {\bf 1}}\int_{\mathbb D^n}\log(1-\overline {a_l}z_l)\prod_{j=1}^n\frac{(1-|a_j|^2)^{k_j}}{\left(1-\overline {a_j}z_j\right)^{2+k_j}}\overline {b(rz)}d\nu(z)\right|.
\Ee
Next we observe that
\Beas
&&\int_{\mathbb D^n}\log(1-\overline {a_l}z_l)\prod_{j=1}^n\frac{(1-|a_j|^2)^{k_j}}{\left(1-\overline {a_j}z_j\right)^{2+k_j}}\overline {b(rz)}d\nu(z)\\ &=& \int_{\mathbb D^n}\log(\frac{1-\overline {a_l}z_l}{1-|a_l|^2})\prod_{j=1}^n\frac{(1-|a_j|^2)^{k_j}}{\left(1-\overline {a_j}z_j\right)^{2+k_j}}\overline {b(rz)}d\nu(z)\\ &+& \int_{\mathbb D^n}\log(1-|a_l|^2)\prod_{j=1}^n\frac{(1-|a_j|^2)^{k_j}}{\left(1-\overline {a_j}z_j\right)^{2+k_j}}\overline {b(rz)}d\nu(z)
\Eeas
and observing that as $\log\left(\frac{1-\overline {a_l}z_l}{1-|a_l|^2}\right)\in \mathbb{B}(\mathbb D^n)$ and $\frac{(1-|a_j|^2)^{k_j}}{\left(1-\overline {a_j}z_j\right)^{2+k_j}}\in A^1(\mathbb D)$ both uniformly, we have by (\ref{eq:testhank1}) that
$$\left|\lim_{r\rightarrow {\bf 1}}\int_{\mathbb D^n}\log\left(\frac{1-\overline {a_l}z_l}{1-|a_l|^2}\right)\prod_{j=1, }^n\frac{(1-|a_j|^2)^{k_j}}{\left(1-\overline {a_j}z_j\right)^{2+k_j}}\overline {b(rz)}d\nu(z)\right|\le C,$$
$l=1,\ldots, n$.
It follows from the latter and (\ref{eq:step1}) that
\Beas
\left(\sum_{j=1}^n\left|\log(1-|a_j|^2)\right|\right)\left(\prod_{j=1}^n(1-|a_j|^2)\right)\left|\lim_{r\rightarrow {\bf 1}}\int_{\mathbb D^n}\prod_{j=1}^n\frac{(1}{\left(1-\overline {a_j}z_j\right)^{2+k_j}}\overline {b(rz)}d\nu(z)\right|\le C<\infty.
\Eeas
That is for any $K\in \N^n$, $$\sup_{z\in \mathbb D^n}\left(\sum_{j=1}^n\left(\log\frac{2}{1-|z_j|^2}\right)\right)\left(\prod_{j=1}^n(1-|z_j|^2)^{k_j}\right)|D^Kb(z)|<\infty.$$
The proof is complete.
\end{proof}

\section{Pointwise multipliers of Bloch spaces}
We prove Theorem \ref{theo:main1} and Theorem \ref{theo:main2} in this section.
\
\begin{proof}[Proof of Theorem \ref{theo:main1}]
We would like to prove that  $\phi$ is such that there is a constant $C>0$ so that for any $f\in \mathbb B(\mathbb D^n)$ and any $z=(z_1,\cdots,z_n)\in \mathbb D^n$,
$$\left[\prod_{j=1}^n(1-|z_j|^2)\right]|D\left(\phi f\right)(z)|<C$$
if and only if $\phi\in H^\infty(\mathbb D^n)\cap \mathbb B_L(\mathbb D^n)$.

We observe that
\Be\label{eq:proddecomposition} D(\phi f)=\phi Df+fD\phi+\sum_{\emptyset \neq J\subset \{1,\ldots,n\},  K={^cJ}}D_{\vec j}\phi D_{\vec k}f\Ee
where $\vec j=(j_1,\ldots,j_l)$ is the vector associated to the set $J=\{j_1,\ldots,j_l\}$ and the same for $\vec k$ and the set $K$ the complement of $J$ in $\{1,\ldots, n\}$.
We also observe that if $J\cap K=\emptyset$, $J\cup K=\{1,\ldots,n\}$ with $J,K\neq \{1,\ldots,n\}$, we have using (\ref{eq:smallblochpartialBloch}) that
\Beas
M &:=& \left[\prod_{j=1}^n(1-|z_j|^2)\right]|D_{\vec j}\phi (z)D_{\vec k}f(z)|\\ &\le& \left(\sup_{z\in \mathbb D^n}\left[\prod_{j\in J}(1-|z_j|^2)\right]|D_{\vec j}\phi(z)|\right)\left(\sup_{z\in \mathbb D^n}\left[\prod_{k\in K}(1-|z_k|^2)\right]|D_{\vec k}f(z)|\right)\\ &\le& \|\phi\|_{\mathbb B}\|f\|_{\mathbb B}.
\Eeas

Remark that if $\phi\in H^\infty(\mathbb D^n)$, then for any $z\in \mathbb D^n$,
$$\left[\prod_{j=1}^n(1-|z_j|^2)\right]|\phi (z)Df(z)|\le \|\phi\|_{\infty}\|f\|_{\mathbb B(\mathbb D^n)}.$$

All the above observations amount to saying that a bounded function $\phi$ is a multiplier from $\mathbb B(\mathbb D^n)$ to $\mathcal B(\mathbb D^n)$ if and only if there is a constant $C>0$ such that for any $f\in \mathbb B(\mathbb D^n)$ and any $z\in \mathbb D^n$,
\Be\label{eq:reductsmalltoprodcond}
\left[\prod_{j=1}^n(1-|z_j|^2)\right]|f(z)D\phi (z)|\le C.
\Ee
Let us suppose that $\phi\in \mathbb B_L(\mathbb D^n)$ and prove that in this case, (\ref{eq:reductsmalltoprodcond}) holds. Using the pointwise estimate of functions in $\mathbb B(\mathbb D^n)$ given by Lemma \ref{lem:pointwisesmallbloch}, we obtain
\Beas \left[\prod_{j=1}^n(1-|z_j|^2)\right]|f(z)D\phi (z)| &\le& \|f\|_{\mathbb B}\left[\prod_{j=1}^n(1-|z_j|^2)\right]\left[\sum_{j=1}^n\log\frac{2}{1-|z_j|^2}\right]|D\phi (z)|\\ &\le& \|f\|_{\mathbb B}\|\phi\|_{\mathbb B_L}.
\Eeas

Conversely, if (\ref{eq:reductsmalltoprodcond}) holds, then testing with the function $f(z)=f_a(z)=\sum_{j=1}^n\log(1-z_j\overline {a_j})$ with $a\in \mathbb D^n$ fixed, we obtain
$$\left[\prod_{j=1}^n(1-|z_j|^2)\right]\left|\sum_{j=1}^n\log(1-z_j\overline {a_j})\right||D\phi (z)|\le C.$$
Taking $z_j=a_j$, $j=1,\ldots, n$ in the latter inequality, we obtain that for any $z\in \mathbb D^n$,
$$\left[\prod_{j=1}^n(1-|z_k|^2)\right]\left[\sum_{j=1}^n\log\frac{2}{1-|z_j|^2}\right]|D\phi (z)|<C$$
which proves that $\phi\in \mathbb B_L(\mathbb D^n)$

\end{proof}

Let us now characterize the multiplier algebra of the product Bloch space. We start with the following elementary result.
\blem\label{lem:derivativeestimate}
Let $\vec j=(j_1,\ldots,j_l)$, $1\le l\le n$, and the associated set $J=\{j_1,\ldots,j_l\}\subseteq \{1,\ldots,n\}$. Then there is a constant $C>0$ such that for any $f\in \mathcal B(\mathbb D^n)$ and any $z=(z_1,\ldots,z_n)\in \mathbb D^n$,
$$\left[\prod_{j=1}^n(1-|z_j|^2)\right]|D_{\vec j}f(z)|\le C\left[\prod_{j\notin J}(1-|z_j|^2)\log\frac{2}{1-|z_j|^2}\right]\|f\|_{\mathcal B}.$$
\elem
\begin{proof}
The proof uses the representation formula of $f\in \mathcal B(\mathbb D^n)$ in Lemma \ref{lem:equivblochprod1}, and \cite[Proposition 1. 4. 10]{Rudin}. We obtain
\Beas
|D_{\vec j}f(z)| &=& \left|\int_{\mathbb D^n}\left[\prod_{j\in J}\frac{2}{(1-z_j\overline {w_j})^3}\right]\left[\prod_{j\notin J}\frac{1}{(1-z_j\overline {w_j})^2}\right]g(w)d\nu(w)\right|\\ &\le& C\|g\|_{\infty}\left\|\left[\prod_{j\in J}\frac{1}{(1-z_j\overline {w_j})^3}\right]\left[\prod_{j\notin J}\frac{1}{(1-z_j\overline {w_j})^2}\right]\right\|_1\\ &\le& C\left[\prod_{j\in J}(1-|z_j|^2)^{-1}\right]\left[\prod_{j\notin J}\log\frac{2}{(1-|z_j|^2)}\right]\|f\|_{\mathcal B}.
\Eeas
The proof is complete.
\end{proof}
As a consequence of the above lemma, we have the following result.
\bprop\label{prop:blochlogblochderivative}
Let $\vec j$ and $\vec k$ be two vectors such that their respective associated sets $J$ and $K$ are complementary in $\{1,\ldots,n\}$, with none of them empty. Suppose that $f\in \mathcal B(\mathbb D^n)$ and $\phi\in \mathcal B_{LL}(\mathbb D^n)$. Then for any $z\in \mathbb D^n$,
$$\left[\prod_{j=1}^n(1-|z_j|^2)\right]|D_{\vec j}f(z)||D_{\vec k}\phi (z)|\lesssim \|f\|_{\mathcal B}\|\phi\|_{\mathcal B_{LL}}.$$
\eprop
\begin{proof}
For $z\in \mathbb D^n$, following Lemma \ref{lem:derivativeestimate} and the definition of $\mathcal B_{LL}(\mathbb D^n)$ we obtain
\Beas
M &:=& \left[\prod_{j=1}^n(1-|z_j|^2)\right]|D_{\vec j}f(z)||D_{\vec k}\phi (z)|\\ &\le& C\|f\|_{\mathcal B}\left[\prod_{j\in K}(1-|z_j|^2)\log\frac{2}{1-|z_j|^2}\right]|D_{\vec k}\phi (z)|\\ &\le& C\|f\|_{\mathcal B}\|\phi\|_{\mathcal B_{LL}}.
\Eeas
\end{proof}
\begin{proof}[Proof of Theorem \ref{theo:main2}]:
\vskip .1cm
Let us start by proving the sufficiency.
Let $\phi$ be as in Theorem \ref{theo:main2}. We would like to prove that for any $f\in \mathcal B(\mathbb D^n)$, the analytic function $\phi f$ belongs to $\mathcal B(\mathbb D^n)$. From the formula (\ref{eq:proddecomposition}), one sees that we only have to estimate the following three terms.
$$M_1:=\left[\prod_{j=1}^n(1-|z_j|^2)\right]|\phi(z)Df(z)|;$$
$$M_2:=\left[\prod_{j=1}^n(1-|z_j|^2)\right]|D_{\vec j}\phi(z) D_{\vec k}f(z)|,$$
$\emptyset \neq J,K\subset \{1,\ldots,n\}=J\cup K,\,\,\, J\cap K=\emptyset;$
$$M_3:=\left[\prod_{j=1}^n(1-|z_j|^2)\right]|f(z)D\phi(z)|.$$
From that $\phi\in H^\infty(\mathbb D^n)$, we obtain
$$M_1:=\left[\prod_{j=1}^n(1-|z_j|^2)\right]|\phi(z)Df(z)|\le \|\phi\|_{\infty}\|f\|_{\mathcal B}.$$
The estimate of the term $M_2$ is given by Proposition \ref{prop:blochlogblochderivative}. For the last term, we use the pointwise estimate of $f\in \mathcal B(\mathbb D^n)$ to get
\Beas M_3 &:=& \left[\prod_{j=1}^n(1-|z_j|^2)\right]|f(z)D\phi(z)|\\ &\le& C\|f\|_{\mathcal B}\left[\prod_{j=1}^n(1-|z_j|^2)\log\frac{2}{1-|z_j|^2}\right]|D\phi(z)|\\ &\le& C\|f\|_{\mathcal B}\|\phi\|_{\mathcal B_L}.
\Eeas
We conclude that for $\phi\in H^\infty(\mathbb D^n)\cap \mathcal B_{LL}(\mathbb D^n)$,
$$\|\phi f\|_{\mathcal B}\le C\left(\|\phi\|_{\infty}+\|\phi\|_{\mathcal B_{LL}}\right)\|f\|_{\mathcal B}.$$
We now prove the necessity part in Theorem \ref{theo:main2}.
We suppose that $\phi$ is such that for any $f\in \mathcal B(\mathbb D^n)$, the function $\phi f$ belongs to $\mathcal B(\mathbb D^n)$. That is there exists a constant $C>0$ such that for any $z\in \mathbb D^n$,
\Be\label{eq:prodinbloch}
(1-\|z\|^2)|D\left(\phi f\right)(z)|=(1-\|z\|^2)\left|\sum_{J\subseteq \{1,\ldots, n\}, K= {^cJ}}D_{\vec j}\phi (z)D_{\vec k}f(z)\right|\le C,
\Ee
with $(1-\|z\|^2)=(1-|z_1|^2)\ldots (1-|z_n|^2)$ and for $J=\emptyset$, $D_{\vec j}f=f$.

We first prove that this necessarily implies that $\phi$ is bounded. Clearly, if $\phi\in \mathcal {M}(\mathcal B(\mathbb D^n))$, then the pointwise estimate of functions in $\mathcal B(\mathbb D^n)$ gives that there exists a constant $C>0$ such that for any $f\in \mathcal B(\mathbb D^n)$ and $z\in \mathbb D^n$,
\Be\label{eq:multipointwise}
\left|f(z)\phi(z)\right|\le C\left(\prod_{j=1}^n\log\frac{2}{1-|z_j|^2}\right)\|f\phi\|_{\mathcal B}\le C\left(\prod_{j=1}^n\log\frac{2}{1-|z_j|^2}\right)\|f\|_{\mathcal B}.
\Ee
We test (\ref{eq:multipointwise}) with the function $$f(z)=f_a(z)=\prod_{j=1}^n\log(1-\overline {a_j}z_j),$$
$a=(a_j,\ldots, a_n)$ given in $\mathbb D^n$. We obtain that for any $z\in \mathbb D^n$,
$$\left|\left(\prod_{j=1}^n\log(1-\overline {a_j}z_j)\right)\phi(z)\right|\le C\left(\prod_{j=1}^n\log\frac{2}{1-|z_j|^2}\right)\|f\|_{\mathcal B}.$$
Taking in particular $z_j=a_j$ ($j=1,\ldots, n$) in the above inequality, we obtain that for any $z\in \mathbb D^n$,
$$|\phi(z)|\le C<\infty,$$
that is $\phi\in H^\infty(\mathbb D^n)$.
\vskip .2cm
We next prove that $\phi\in \mathcal B_{LL}(\mathbb D^n)$. For this, we first observe the following fact.
\blem\label{lem:loweringthedim}
If $\phi\in \mathcal {H}(\mathbb D^n)$ is a multiplier of $\mathcal B(\mathbb D^n)$, then for any fixed $a_1\in \mathbb D$, the function $\phi(a_1,\cdot)$ is a multiplier of $\mathcal B(\mathbb D^{n-1})$. Moreover, $$\|\phi(a_1,\cdot)|_{\mathcal B(\mathbb D^{n-1})\rightarrow \mathcal B(\mathbb D^{n-1})}\lesssim \|\phi|_{\mathcal B(\mathbb D^{n})\rightarrow \mathcal B(\mathbb D^{n})}.$$
\elem
\begin{proof}
We first prove that for any $w_n\in \mathbb D$ fixed, for any $b\in \mathcal B(\mathbb D^n)$, the function $b(\cdot,w_n)$ which is a function of $n-1$ variables, is in $\mathcal B(\mathbb D^{n-1})$ with
\Be\label{eq:prodblochinonepara}\|b(\cdot,w_n)\|_{\mathcal B(\mathbb D^{n-1})}\lesssim \log\frac{4}{1-|w_n|^2}\|b\|_{\mathcal B(\mathbb D^n)}.
\Ee
Let $z=(z_1,\ldots,z_{n-1})\in \mathbb D^{n-1}$. From the integral representation of elements of $\mathcal B(\mathbb D^n)$, we have that for some $g\in L^\infty(\mathbb D^n)$,
$$b(z,w_n)=\int_{\mathbb D^n}\frac{g(\xi)d\nu(\xi)}{(1-z_1\bar {\xi_1})^2(1-z_2\bar {\xi_2})^2\ldots (1-z_{n-1}\bar {\xi_n})^2(1-w_n\bar {\xi_n})^2},$$
hence
\Beas
&&\left(\prod_{j=1}^{n-1}(1-|z_j|^2)\right)\left|D_1\ldots D_{n-1}b(z,w_n)\right|\\ &=& \left(\prod_{j=1}^{n-1}(1-|z_j|^2)\right)\left|\int_{\mathbb D^n}\frac{2^ng(\xi)d\nu(\xi)}{\left(\prod_{j=1}^{n-1}(1-z_j\bar {\xi_j})^3\right)(1-w_n\bar {\xi_n})^2}\right|\\ &\le& 2^n\left(\prod_{j=1}^{n-1}(1-|z_j|^2)\right)\int_{\mathbb D^n}\frac{|g(\xi)|d\nu(\xi)}{\left(\prod_{j=1}^{n-1}|1-z_j\bar {\xi_j}|^3\right)|1-w_n\bar {\xi_n}|^2}\\ &\lesssim& \left(\prod_{j=1}^{n-1}(1-|z_j|^2)\right)\|g\|_{L^\infty(\mathbb D^n)}\int_{\mathbb D^n}\frac{d\nu(\xi)}{\left(\prod_{j=1}^{n-1}|1-z_j\bar {\xi_j}|^3\right)|1-w_n\bar {\xi_n}|^2}\\ &\lesssim& \|b\|_{\mathcal B(\mathbb D^n)}\log\frac{4}{1-|w_n|^2}.
\Eeas
Now let $\phi$ be a multiplier of $\mathcal B(\mathbb D^n)$. Then from (\ref{eq:prodblochinonepara}) we obtain that for any $b\in \mathcal B(\mathbb D^n)$ and any $w_n\in \mathbb D$ fixed,
$$\left(\prod_{j=1}^{n-1}(1-|z_j|^2)\right)\left|D_1\ldots D_{n-1}(\phi b)(z,w_n)\right| \lesssim \left(\log\frac{4}{1-|w_n|^2}\right)\|\phi b\|_{\mathcal B(\mathbb D^n)} $$
and so
\Be\label{eq:prodblochinoepara1}
\left(\prod_{j=1}^{n-1}(1-|z_j|^2)\right)\left|D_1\ldots D_{n-1}(\phi b)(z,w_n)\right| \lesssim \left(\log\frac{4}{1-|w_n|^2}\right)\|\phi\|_{\mathcal B(\mathbb D^n)\rightarrow \mathcal B(\mathbb D^n)}\|b\|_{\mathcal B(\mathbb D^n)}
\Ee

Let us take in (\ref{eq:prodblochinoepara1}), $b(z,\xi)=g(z)\log(1-\xi\bar {w_n})$, $g\in \mathcal {B}(\mathbb D^{n-1})$, $z\in \mathbb D^{n-1}$ and $\xi\in \mathbb D$. We obtain
\Beas
S &:=& \left(\prod_{j=1}^{n-1}(1-|z_j|^2)\right)|\log(1-\xi\bar {w_n})|\left|D_1\ldots D_{n-1}(\phi g)(z,w_n)\right|\\ &=& \left(\prod_{j=1}^{n-1}(1-|z_j|^2)\right)\left|D_1\ldots D_{n-1}(\phi b)(z,w_n)\right|\\ &\lesssim&  \log\frac{4}{1-|w_n|^2}\|\phi\|_{\mathcal B(\mathbb D^n)\rightarrow \mathcal B(\mathbb D^n)}\|g\|_{\mathcal B(\mathbb D^{n-1})}.
\Eeas
Taking $\xi=w_n$ in the above inequalities, we obtain that for any $g\in \mathcal {B}(\mathbb D^{n-1})$, and any $z\in \mathbb D^{n-1}$,
$$\left(\prod_{j=1}^{n-1}(1-|z_j|^2)\right)\left|D_1\ldots D_{n-1}(\phi g)\right|\lesssim \|M_\phi\|_{\mathcal B(\mathbb D^n)\rightarrow \mathcal B(\mathbb D^n)}\|g\|_{\mathcal B(\mathbb D^{n-1})}.$$
Thus for any $w_n\in \mathbb D$ fixed, $\phi(\cdot,w_n)$ is a multiplier of $\mathcal B(\mathbb D^{n-1})$. The proof of the lemma is complete.
\end{proof}
We next proceed by induction on the number of parameters $n\ge 2$ to prove that if $\phi$ is a multiplier of $\mathcal B(\mathbb D^n)$, then $\phi\in \mathcal {B}_{LL}(\mathbb D^n)$. We start by the case $n=2$. Let $\phi$ be a multiplier of $\mathcal B(\mathbb D^2)$. Then there exists a constant $C>0$ such that for any $b\in \mathcal B(\mathbb D^2)$ and any $z=(z_1,z_2)$,
\Be\label{eq:first}(1-|z_1|^2)(1-|z_2|^2)\left|D\phi f(z)\right|\le C\|f\|_{\mathcal B(\mathbb D^2)}.
\Ee
Recall that in this case,  $D=D_1D_2$ and $$D(\phi f)(z)=f(z)D\phi(z)+\phi(z)Df(z)+D_1\phi(z)D_2(z)+D_1f(z)D_2\phi(z).$$
But by Lemma \ref{lem:loweringthedim} and Proposition \ref{prop:multionepara}, $\phi(\cdot,z_2)$ and $\phi(z_1,\cdot)$ are uniformly in $\mathcal {B}_L(\mathbb D)$, that is there is a constant $C>0$ such that for any $z=(z_1,z_2)\in \mathbb D^2$,
$$(1-|z_1|^2)\left(\log\frac{4}{1-|z_1|^2}\right)|D_1\phi (z)|\le C$$
and $$(1-|z_2|^2)\left(\log\frac{4}{1-|z_2|^2}\right)|D_2\phi(z)|\le C.$$
Hence for any $z=(z_1,z_2)\in \mathbb D^2$, we obtain using Lemma \ref{lem:derivativeestimate},
\Beas
S &:=& (1-|z_1|^2)(1-|z_2|^2)|D_1\phi(z)||D_2f(z)|\\ &\lesssim& (1-|z_1|^2)\left(\log\frac{4}{1-|z_1|^2}\right)|D_1\phi(z)|\|f\|_{\mathcal B(\mathbb D^2)}
\Eeas
and consequently,
\Be\label{eq:first1}
(1-|z_1|^2)(1-|z_2|^2)|D_1\phi(z)||D_2f(z)| \le C\|f\|_{\mathcal B(\mathbb D^2)}.
\Ee
In the same way, we obtain for any $z=(z_1,z_2)\in \mathbb D^2$,
\Be\label{eq:first2}(1-|z_1|^2)(1-|z_2|^2)|D_2\phi(z)||D_1f(z)|\le C\|f\|_{\mathcal B(\mathbb D^2)}.
\Ee
Also, note that as $\phi\in H^\infty(\mathbb D^2)$, we have that for any $z=(z_1,z_2)\in \mathbb D^2$,
\Be\label{eq:first3}(1-|z_1|^2)(1-|z_2|^2)|\phi(z)||Df(z)|\le \|\phi\|_{\infty}\|f\|_{\mathcal B(\mathbb D^2)}.
\Ee
From (\ref{eq:first}), (\ref{eq:first1}), (\ref{eq:first2}) and (\ref{eq:first3}), we deduce that there exists a constant $C>0$ such that for any $f\in \mathcal B(\mathbb D^2)$ and any $z=(z_1,z_2)\in \mathbb D^2$,
\Be\label{eq:first4}
(1-|z_1|^2)(1-|z_2|^2)|f(z)||D\phi(z)|\le C\|f\|_{\mathcal B(\mathbb D^2)}.
\Ee
For $a=(a_1,a_2)\in \mathbb D^2$ given, we test (\ref{eq:first4}) with $$f(z)=f_a(z)=\log(1-z_1\bar {a_1})\log(1-z_2\bar {a_2})$$ which is uniformly in  $\mathcal B(\mathbb D^2)$ and obtain for any $z=(z_1,z_2)\in \mathbb D^2$,
\Be\label{eq:first5}
(1-|z_1|^2)(1-|z_2|^2)|\log(1-z_1\bar {a_1})||\log(1-z_2\bar {a_2})||D\phi(z)|\le C\|f\|_{\mathcal B(\mathbb D^2)}.
\Ee
Taking in particular $z_1=a_1$ and $z_2=a_2$ in (\ref{eq:first5}), we conclude that there is a constant $C>0$ such that for any $z=(z_1,z_2)\in \mathbb D^2$,
$$(1-|z_1|^2)(1-|z_2|^2)\left(\log\frac{4}{1-|z_1|^2}\right)\left(\log\frac{4}{1-|z_2|^2}\right)|D\phi(z)|\le C,$$
that is $\phi\in \mathcal B_L(\mathbb D^2)$. This completes the proof for the case $n=2$.
\vskip .2cm
Now for $n>2$, we suppose that $\phi$ is a multiplier of $\mathcal B(\mathbb D^n)$ implies that $\phi \in \mathcal B_{LL}(\mathbb D^n)$. We prove that this implies that if $\phi$ is a multiplier of $\mathcal B(\mathbb D^{n+1})$, then $\phi\in \mathcal B_{LL}(\mathbb D^{n+1})$.
\vskip .1cm
Let $\phi$ be a multiplier of  $\mathcal B(\mathbb D^{n+1})$. Then by Lemma \ref{lem:loweringthedim}, for any $w_{n+1}\in \mathbb D$ fixed, $\phi(\cdot,w_{n+1})$ is a multiplier of $\mathcal B(\mathbb D^{n})$ with uniformly bounded multiplier norm. Hence by our hypothesis, $\phi(\cdot,w_{n+1})\in \mathcal B_{LL}(\mathbb D^{n})$ uniformly. It follows in particular that there is a constant $C>0$ such that for any $\vec j=(j_1,\ldots,j_l)$ with associated set $J=\{j_1,\ldots,j_l\}\subset \{1,2,\ldots,n+1\}$, and any $z=(z_1,\ldots,z_{n+1})\in \mathbb D^{n+1}$,
\Be\label{eq:uniformlowering}\left(\prod_{j\in J }(1-|z_j|^2)\log\frac{4}{1-|z_j|^2}\right)|D_{\vec j}\phi(z)|\le C.
\Ee
 Denoting by $K$ the complement set of $J$ in $\{1,2,\ldots,n+1\}$ with associated vector $\vec k$, we obtain using Lemma \ref{lem:derivativeestimate} that for any $f\in \mathcal B(\mathbb D^{n+1})$, and any $z=(z_1,\ldots,z_{n+1})\in \mathbb D^{n+1}$,
\Beas\label{eq:endeq}
Q &:=&\left(\prod_{j=1}^{n+1}(1-|z_j|^2)\right)|D_{\vec j}\phi(z)||D_{\vec k}f(z)|\\ &\le& C\|f\|_{\mathcal B(\mathbb D^{n+1})}\left(\prod_{j\in J }(1-|z_j|^2)\log\frac{4}{1-|z_1|^2}\right)|D_{\vec j}\phi(z)|.
\Eeas
Hence applying (\ref{eq:uniformlowering}) to the above, we obtain
\Be\label{eq:control1}
\left(\prod_{j=1}^{n+1}(1-|z_j|^2)\right)|D_{\vec j}\phi(z)||D_{\vec k}f(z)|\le  C\|f\|_{\mathcal B(\mathbb D^{n+1})}.
\Ee
Also we have since $\phi\in H^\infty(\mathbb D^{n+1})$, that for any $z=(z_1,\cdots,z_{n+1})\in \mathbb D^{n+1}$,
\Be\label{eq:endeq1}
\left(\prod_{j=1}^{n+1}(1-|z_j|^2)\right)|\phi(z)||Df(z)|\le \|\phi\|_{H^\infty}\|f\|_{\mathcal B(\mathbb D^{n+1})}.
\Ee
We recall that in this case, $$D(\phi f)=fD\phi+\phi Df+\sum_{\emptyset \neq J\subset \{1,\cdots,n+1\}}D_{\vec j}\phi D_{\vec k}f .$$ From (\ref{eq:control1}), (\ref{eq:endeq1}) and the fact that we have a constant $C>0$ such that for any $f\in \mathcal B(\mathbb D^{n+1})$ and any $z=(z_1,\ldots,z_{n+1})\in \mathbb D^{n+1}$,
\Be
\left(\prod_{j=1}^{n+1}(1-|z_j|^2)\right)\left|D\left(\phi f\right)(z)\right|\le C\|f\|_{\mathcal B(\mathbb D^{n+1})},
\Ee
we obtain that there exists a constant $C>0$ such that for any $f\in \mathcal B(\mathbb D^{n+1})$ and for any $z=(z_1,\ldots,z_{n+1})\in \mathbb D^{n+1}$,
\Be\label{eq:endeq2}
\left(\prod_{j=1}^{n+1}(1-|z_j|^2)\right)|f(z)||D\phi(z)|\le C\|f\|_{\mathcal B(\mathbb D^{n+1})}.
\Ee
For $a=(a_1,a_2,\ldots,a_{n+1})\in \mathbb D^{n+1}$ given, we test (\ref{eq:endeq2}) with $$f(z)=f_a(z)=\log(1-z_1\bar {a_1})\ldots\log(1-z_{n+1}\bar {a_{n+1}})$$ which is uniformly in  $\mathcal B(\mathbb D^{n+1})$ and obtain for any $z=(z_1,\ldots,z_{n+1})\in \mathbb D^{n+1}$,
\Be\label{eq:endeq3}
\left(\prod_{j=1}^{n+1}(1-|z_j|^2)\left|\log(1-z_j\bar {a_j})\right|\right)|D\phi(z)|\le C.
\Ee
Taking in particular $z_j=a_j$, $j=1,2,\ldots,n+1$ (\ref{eq:endeq3}), we obtain that there is a constant $C>0$ such that for any $z=(z_1,\cdots, z_{n+1})\in \mathbb D^{n+1}$,
$$
\left(\prod_{j=1}^{n+1}(1-|z_j|^2)\left(\log\frac{4}{1-|z_j|^2}\right)\right)|D\phi(z)|\le C.
$$
That is $\phi\in \mathcal {B}_L(\mathbb D^{n+1})$. The latter and (\ref{eq:uniformlowering}) allow us to conclude that $\phi\in \mathcal {B}_{LL}(\mathbb D^{n+1})$. The proof is complete.
\end{proof}
\section{Remarks on the pointwise Bloch space}
The multiplier algebra of $\mathbb B(\mathbb D^n)$ has been found by F. Colonna and R. F. Allen in \cite{AlCol}. They proved exactly the following.
\bprop\label{prop:AlCol}
The only multipliers of $\mathbb B(\mathbb D^n)$ are the constants.
\eprop

 We have the following consequence of the above proposition.
\bcor
Let $X$ be a Banach space of analytic functions strictly containing $\mathbb B(\mathbb D^n)$. Then $\mathcal {M}\left(X,\mathbb B(\mathbb D^n)\right)=\{0\}$.
\ecor
In particular, we obtain the following.
\bcor
$\mathcal {M}\left(\mathcal B(\mathbb D^n),\mathbb B(\mathbb D^n)\right)=\{0\}$.
\ecor

\emph{Acknowledgement}  This work was partly completed at the Centre for Advanced Study of Oslo within the research group "Operator related function theory and time-frequency analysis". I would like to thank professors Yurii Lyubarskii and Kristian Seip for their kind invitation and the centre for support.
\vskip 1cm

\section*{Conflict of interest statement}
The author declares that there is no conflict of interests
regarding the publication of this paper.


  \bibliographystyle{elsarticle-num} 

\begin{thebibliography}{00}


\bibitem{AlCol}
\textsc{R. F.  Allen,  F. Colonna},  \emph{Multiplication operators on the Bloch space of bounded homogeneous domains}. Comput. Methods Funct. Theory {\bf 9} (2009), no. 2, 679--693.

\bibitem{BL}
\textsc{A.~Bonami, Luo~ Luo}, \emph{On Hankel operators between Bergman spaces on the unit ball}. Houston J. Math. {\bf 31} (2005), no. 3, 815--827.

\bibitem{Constantin}
\textsc{O. Constantin}, \emph{Weak product decompositions and Hankel operators on vector-valued Bergman spaces}. J. Oper. Theor. {\bf 59} (2008), no. 1, 157–-178.




\bibitem{JPR}
\textsc{S.Janson, J.Peetre, R.Rochberg}, \emph{Hankel forms and the Fock space}. Revista Math. Ibero-Amer. {\bf 3} (1987) 61–138.

\bibitem{Har}
\textsc{A. Harutyunyan}, \emph{Bloch spaces of holomorphic functions in the polydisk}. J. Funct. Spaces Appl. {\bf 5} (2007), no. 3,  213–-230.

\bibitem{HarLus}
\textsc{A. Harutyunyan, W. Lusky}, \emph{Weighted holomorphic Besov spaces on the polydisk}. J. Funct. Spaces Appl. {\bf 9} (2011), no. 1, 1--16.

\bibitem{Rudin}
\textsc{W.~Rudin}, \emph{Function theory in the unit ball of
$\C^n$}. Grundlehren der Mathematischen Wissenschaften [Fundamental Principles of Mathematical Science], 241.
Springer-Verlag, New York-Berlin (1980).



\bibitem{benoit1}
\textsc{B.~F.~Sehba}, \emph{On some equivalent definitions of $\rho$-Carleson measures on the unit ball}. Acta Sci. Math. (Szeged) {\bf 75} (2009), no. 3--4, 499--525.

\bibitem{stevo}
\textsc{S.~Stevic}, \emph{A note on a Theorem of Zhu on weighted Bergman projections on the polydisc}. Houston J. Math. {\bf 34} (2004), no. 4, 511--521.

\bibitem{Timoney}
 \textsc{R. M. Timoney},  \emph{Bloch functions in several complex variables I}. Bull.  London Math. Soc. {\bf 12}
 (1980), 241--267.

\bibitem{zhao} \textsc{R.~Zhao}, \emph{On logarithmic Carleson measures}, Acta Sci.~Math (Szeged) {\bf 69} (2003), no 3-4, 605--618
\bibitem{Zhu1}
\textsc{K. Zhu,} \emph{Weighted Bergman projections on the polydisc}. Houston J. Math. {\bf 20} (1994), no. 2, 275–-292.

\bibitem{Zhu}
\textsc{K. Zhu,} \emph{Spaces of holomorphic functions in the unit ball}. Graduate Texts in Mathematics 226,Springer  Verlag (2004).

 \end{thebibliography}


\end{document}